\documentclass[12pt,english]{article}
\usepackage{amsfonts,here,graphicx}
\usepackage{cite,amsfonts,amssymb} 
\usepackage[notref,notcite]{showkeys}

\def\proof{\noindent \medskip {\bf Proof:}$\;\;$}

\def\ass#1#2\endass{\vskip5pt plus2pt \noindent{\bf (A.#1)} {\it #2} \vskip5pt
plus2pt }

\newtheorem{prop}{Proposition}

\newtheorem{defin}{Definition}

\title{\bf Stochastic Integral with respect to Cylindrical Wiener 
Process}
\author{\sf Anna Karczewska \\ \\ \it
Institute of Mathematics, Maria Curie--Sk{\l}odowska University\\ \it
pl. M. Curie--Sk{\l}odowskiej 1, PL--20--031 Lublin, Poland}
\date{}

\begin{document}

\maketitle

 \noindent\def\thefootnote{}
\footnotetext{
{\em 1991 Mathematics Subject Classification:}
Primary: 60H05; Secondary: 60H30.\\
{\em Key words and phrases: Stochastic integral, 
infinite dimensional Wiener process, cylindrical Wiener process} }

\hfill{\it Dedicated to Professor Dominik Szynal
 on the occasion ofhis 60--th birthday}\\

\begin{abstract}
This paper is devoted to a construction of the stochastic It\^o integral
with respect to infinite dimensional cylindrical Wiener process. The
construction given is an alternative one to that introduced by DaPrato and
Zabczyk \cite{3}. The connection of the introduced integral with the
integral defined by Walsh \cite{9} is provided as well.
\end{abstract}

\section{Introduction}

Recently there have been written several papers devoted to stochastic
partial differential equations forced by cylindrical Wiener process, e.g.,
\cite{4}, \cite{2} and \cite{7}. In the study of stochastic partial
diffrential equations some authors (see references given in Chapter
4 in \cite{3}) have used a stochastic integral with respect 
to the so--called Brownian sheet, which is a special kind of
 cylindrical Wiener
process, rather than with respect to cylindrical Wiener 
process in general form.

In the paper we provide a construction of stochastic integral
with respect to an infinite dimensional cylindrical Wiener process
alternative to the
construction given by DaPrato and Zabczyk in their monograph \cite{3}.
We introduce the convenient construction which is based, by analogy
to the construction given by Ichikawa \cite{6} for the integral
with respect to classical infinite dimensional Wiener process, on the
stochastic integrals with respect to real--valued Wiener processes.
The advantage
of using of such a construction is that we can use basic results and
arguments of the  finite dimensional case.
Finally, we compare the integral constructed in the paper with the
integral introduced by Walsh \cite{9}.

Let us recall from \cite{6} the definition of
Wiener process with values in Hilbert space $U$
(called later the classical infinite dimensional Wiener process)
and the stochastic integral
with respect to this Wiener process.

\begin{defin} \label{def1} Let $Q:\ U \to U$ be a linear symmetric
non--negative nuclear operator ($\mathrm{tr} Q<+\infty$). 
A square integrable
$U$--valued stochastic process $W(t)$, $t\ge0$,
defined on a probability space $(\Omega,\mathcal{F},
(\mathcal{F}_t)_{t\ge0},P)$,
where $\mathcal{F}_t$ denote $\sigma$--fields such that
$\mathcal{F}_t\subset \mathcal{F}_s \subset \mathcal{F}$ for $t<s$, 
is called {\tt Wiener process} with covariance operator $Q$ if:
\begin{enumerate}
\item $W(0)=0$,
\item $EW(t)=0$, $\mathrm{Cov} [W(t)-W(s)]=(t-s)Q$ for all $s,t\ge 0$,
\item $W$ has independent increments,
\item $W$ has continuous trajectories,
\item $W$ is adapted with respect to the filtraction $(\mathcal{F}_t)$, that
is, for any $t\ge0$, $W(t)$ is $\mathcal{F}_t$--measurable.
\end{enumerate}
\end{defin}

In the light of the above, Wiener process is Gaussian and has the
following structure:
let $\{d_i\}\subset U$ be an orthonormal set of eigenvectors of $Q$ with
corresponding eigenvalues $\zeta_i$ (so $\mathrm{tr} Q=\sum_{i=1}^{\infty}\zeta_i$),
then $W(t)=\sum_{i=1}^{\infty}\beta_i(t)d_i$, where $\beta_i$ are independent
real Wiener processes with $E(\beta_i^2(t))=\zeta_it$.
This type of structure of Wiener process will be used in definition of the
stochastic integral.
\vskip2mm
Let $L(U,Y)$ denote the space of linear bounded operators from $U$ into
$Y$.
\vskip2mm
For any Hilbert space $Y$ we denote by $M(Y)$ the space of all stochastic
processes
$
g:\ [0,T]\times \Omega \rightarrow L(U,Y)
$
such that
$$
E\left(\int_0^T\|g(t)\|_{L(U,Y)}^2dt\right)<+\infty
$$
and for all $u\in U$, $g(t)u$ is a $Y$--valued stochastic process measurable
with respect to the filtration $(\mathcal{F}_t)$.
\vskip2mm
The stochastic integral $\int_0^tg(s)dW(s)\in Y$ is defined for all
$g\in M(Y)$ by
$$
\int_0^tg(s)dW(s)=\lim_{m\to\infty}\sum_{i=1}^m\int_0^tg(s)d_id\beta_i(s)
$$
in $L^2(\Omega)$ sense.

We shall show that the series in the above formula is
convergent.
\vskip5mm
Let
$W^{(m)}(t)=\sum_{i=1}^md_i\beta_i(t).$
Then, the integral
$$
\int_0^tg(s)dW^{(m)}(s)=\sum_{i=1}^m\int_0^tg(s)d_id\beta_i(s)
$$
is well defined for $g\in M(Y)$ and additionally
$$
\int_0^tg(s)dW^{(m)}(s) ~{\longrightarrow}_{{\hspace{-4ex}}_{m\to\infty}}
\int_0^tg(s)dW(s)\hskip5mm \mbox{in}\ Y
$$
in $L^2(\Omega)$ sense.

This convergence comes from the fact that the sequence
$$
y_m=\int_0^tg(s)dW^{(m)}(s),\qquad m\in \mathbb{N}
$$
is Cauchy sequence in the space of square integrable random variables.
For, using features of stochastic integrals with respect to
$\beta_i(s)$, for any $m,n\in \mathbb{N}$, $m<n$, we have:
\begin{eqnarray}
E\left(\|y_n-y_m\|_Y^2\right) &=&
\sum_{i=m+1}^n\zeta_iE\int_0^t\left(g(s)d_i,g(s)d_i\right)_Yds \\
&\leq& \left(\sum_{i=m+1}^n\zeta_i\right)E\int_0^t\|g(s)\|_{L(U,Y)}^2ds
 ~~\longrightarrow_{{\hspace{-4ex}}_{m,n\to\infty}}~ 0. \nonumber
\end{eqnarray}
Hence, there exists a limit of the sequence $(y_m)$ which defines the stochastic
integral $\int_0^tg(s)dW(s)$.
\vskip2mm
The above construction of the stochastic integral required the assumption that
$Q$ was a nuclear operator. (This assumption was used in (1).) However, it
is possible to extend the definition of the stochastic integral to the case
of general bounded self-adjoint, non--negative operator $Q$ on Hilbert space
$U$. (But it will require some restrictions on the integrand $g$.) Stochastic
integral for this case has been defined e.g. in the monograph \cite{3}.
(To avoid trivial complications we shall assume that $Q$ is strictly
positive, that is: $Q$ is non--negative and $Qx\ne 0$ for $x\ne 0$.)
\vskip2mm
Let us recall the following definition.
\vskip5mm
\begin{defin} \label{def2} 
(\cite{1} or \cite{3}) Let $E$ and $F$ be separable Hilbert spaces with
orthonormal bases $\{e_k\}\subset E$ and $\{f_j\}\subset F$,
respectively. A linear bounded operator $T:\ E\to F$ is called
{\tt Hilbert-Schmidt operator} if $\sum_{k=1}^{\infty}\|Te_k\|_F^2
<+\infty$.
\end{defin}
\vskip2mm
Because
$$
\sum_{k=1}^{\infty}\|Te_k\|_F^2=\sum_{k=1}^{\infty}\sum_{j=1}^{\infty}
(Te_k,f_j)_F^2=\sum_{j=1}^{\infty}\|T^*f_j\|_E^2,
$$
where $T^*$ denotes the operator adjoint to $T$, then the definition of
Hilbert-Schmidt operator and the number
$\|T\|_{HS}=\left(\sum_{k=1}^{\infty} \|Te_k\|_F^2\right)^{\frac12}$
do not depend on the basis $\{e_k\}$, $k\in\mathbb{N}$.
Moreover
$\|T\|_{HS}=\|T^*\|_{HS}.$

Additionally, $L_2(E,F)$ -- the set of all Hilbert-Schmidt operators
from $E$ into $F$, endowed with the norm $\|\cdot \|_{HS}$ defined
above, is a separable Hilbert space.
\vskip2mm
Let us introduce the subspace $U_0$ of the space $U$ defined by
$U_0=Q^{\frac12}(U)$ with the norm
$$
\|u\|_{U_0}=\Vert Q^{-\frac12}u\Vert_U,\qquad u\in U_0.
$$

Assume that $U_1$ is an arbitrary Hilbert space such that $U$ is
continuously embedded
into $U_1$ and the embedding of $U_0$ into $U_1$ is a Hilbert-Schmidt operator.
\vskip2mm
In particular
\begin{enumerate}
\item When $Q=I$, then $U_0=U$ and the embedding of $U$ into $U_1$ is
Hilbert-Schmidt operator.
\item When $Q$ is a nuclear operator, that is $\mathrm{tr} Q<+\infty$, then
$U_0=Q^{\frac12}(U)$ and we can take $U_1=U$. Because in this case
$Q^{\frac12}$ is Hilbert-Schmidt operator then the embedding $U_0\subset U$ is
Hilbert-Schmidt operator.
\end{enumerate}

\section{Stochastic integral with respect to cylindrical Wiener process}
 
We denote by $L_2^0=L_2(U_0,Y)$ the space of Hilbert-Schmidt 
operators acting
from $U_0$ into $Y$, and by $L=L(U,Y)$, like earlier, we denote
the space of linear bounded operators
from $U$ into $Y$.
\vskip2mm
Let us consider the norm of the operator $\psi \in L_2^0$:

\begin{eqnarray*}
\Vert \psi\Vert_{L_2^0}^2&=&
\sum_{h,k=1}^{\infty}\left(\psi g_h,f_k\right)_Y^2
=\sum_{h,k=1}^{\infty}\lambda_h\left(\psi e_h,f_k\right)_Y^2\\
&=&\Vert \psi Q^{\frac12}\Vert_{HS}^2=\mathrm{tr} (\psi Q\psi^*),
\end{eqnarray*}
where $g_j=\sqrt{\lambda_j}e_j$, and $\{\lambda_j\}$, $\{e_j\}$ are
eigenvalues and eigenfunctions of the operator $Q$;
\newline
$\{g_j\}$, $\{e_j\}$ and
$\{f_j\}$ are orthonormal bases of spaces $U_0$, $U$ and $Y$,
respectively.
\vskip2mm
The space $L_2^0$ is a separable Hilbert space with the norm
$\Vert\psi\Vert_{L_2^0}^2=\mathrm{tr}\left(\psi Q\psi^*\right)$.

In particular
\begin{enumerate}
\item When $Q=I$ then $U_0=U$ and the space $L_2^0$ becomes $L_2(U,Y)$.
\item When $Q$ is a nuclear operator, that is $\mathrm{tr} Q<+\infty$, then
$L(U,Y)\subset L_2(U_0,Y)$. For,
assume that $K\in L(U,Y)$ that is $K$ is linear
bounded operator from the space $U$ into $Y$. Let us consider the operator
$\psi =K|_{U_0}$, that is the restriction of operator $K$ to the space $U_0$,
where $U_0=Q^{\frac12}(U)$. Because $Q$ is nuclear operator, then $Q^{\frac12}$
is Hilbert-Schmidt operator. So, the embedding $J$ of the space $U_0$ into
$U$ is Hilbert-Schmidt operator. We have to compute the norm
$\Vert\psi\Vert_{L_2^0}$ of the operator $\psi:\ U_0\to Y$. We obtain
$\Vert\psi\Vert_{L_2^0}^2\equiv \Vert KJ\Vert_{L_2^0}^2=
\mathrm{tr} KJ(KJ)^*,$
where $J:\ U_0\to U$.
\vskip2mm
Because $J$ is Hilbert--Schmidt operator and $K$ is linear bounded operator
then, basing on the theory of Hilbert--Schmidt operators (e.g. \cite{5},
Chapter I), $KJ$ is Hilbert--Schmidt operator, too. Next, $(KJ)^*$ is
Hilbert--Schmidt operator. In consequence, $KJ(KJ)^*$ is nuclear operator, so
$\mathrm{tr} KJ(KJ)^*<+\infty$. Hence, $\psi =K|_{U_0}$ is Hilbert-Schmidt operator
on the space $U_0$, that is $K\in L_2(U_0,Y)$.
\end{enumerate}

Let $\{g_j\}$ denote an orthonormal basis in $U_0$ and $\{\beta_j\}$ be
a family of independent standard real-valued Wiener processes.
\vskip2mm
Although Propositions 1. and 2. introduced below are known (see, e.g.
Proposition 4.11 in the monograph \cite{3}), because of their importance
we formulate them again and provide with detailed proofs.
 
\begin{prop} \label{prop1}
The formula
\begin{equation}
W_c(t)=\sum_{j=1}^{\infty}g_j\beta_j(t),\qquad t\ge 0
\end{equation}
defines Wiener process in $U_1$ with covariance operator $Q_1$ such that
$\mathrm{tr} Q_1<+\infty$.
\end{prop}
\proof{
This comes from the fact that the series (2) is convergent in space
$L^2(\Omega,\mathcal{F}, P;U_1)$. We have
\begin{eqnarray*}
&~& E\left(\left\Vert\sum_{j=1}^ng_j\beta_j(t)-\sum_{j=1}^mg_j\beta_j(t)
\right\Vert_{U_1}^2
\right)=E\left(\left\Vert \sum_{j=m+1}^n g_j\beta_j(t)
\right\Vert_{U_1}^2\right) = \\
&=& E\left(\sum_{j=m+1}^ng_j\beta_j(t),\ \sum_{k=m+1}^ng_k\beta_k(t)\right)
_{U_1}
=E\sum_{j=m+1}^n\left(g_j\beta_j(t),g_j\beta_j(t)\right)_{U_1}\\
&=&E\left(\sum_{j=m+1}^n(g_j,g_j)_{U_1}\beta_j^2(t)\right)
=t\sum_{j=m+1}^n\Vert g_j\Vert_{U_1}^2,\qquad n\ge m\ge 1.
\end{eqnarray*}
From the  assumption, the embedding $J:\ U_0\to U_1$ is Hilbert--Schmidt operator,
then for the basis $\{g_j\}$, complete and orthonormal in $U_0$, we have
$\sum_{j=1}^{\infty}\left\Vert Jg_j\right\Vert_{U_1}^2<+\infty.$
Because $Jg_j=g_j$ for any $g_j\in U_0$, then
$\sum_{j=1}^{\infty}\left\Vert g_j\right\Vert_{U_1}^2\!<\!+\!\infty$
which means
$\sum_{j=m+1}^n\left\Vert g_j\right\Vert_{U_1}^2 \to 0$
when $m,n\to \infty.$
\vskip1mm
Conditions 1), 2), 3) and 5) of the definition of Wiener process are
obviously satisfied. The process defined by (2) is Gaussian because
$\beta_j(t),\ j\in\mathbb{N}$, are independent Gaussian processes.
By Kolmogorov test theorem (see, e.g. \cite{3}, Theorem 3.3), trajectories of
the process $W_c(t)$ are continuous (condition 4) of the definition of Wiener
process) because $W_c(t)$ is Gaussian.

Let $Q_1:\ U_1\to U_1$ denote the covariance operator of the process $W_c(t)$
defined by (2). From the definition of covariance, for $a,b\in U_1$ we have:
\begin{eqnarray*}
\left( Q_1 a,b \right)_{_{U_1}}&=&E(a,W_c(t))_{_{U_1}}(b,W_c(t))_{_{U_1}}
=E\left(\sum_{j=1}^{\infty}(a,g_j)_{_{U_1}}(b,g_j)_{_{U_1}}
\beta_j^2(t)\right)\\
&=&t\sum_{j=1}^{\infty}(a,g_j)_{_{U_1}}(b,g_j)_{_{U_1}}
=t\left(\sum_{j=1}^{\infty}g_j(a,g_j)_{_{U_1}},b\right)_{U_1}.
\end{eqnarray*}
Hence $Q_1a=t\sum_{j=1}^{\infty}g_j(a,g_j)_{_{U_1}}.$
\vskip2mm
Because the covariance operator $Q_1$ is non--negative, then
(by Proposition C.3 in \cite{3}) $Q_1$ is a nuclear operator if and only if
$\sum_{j=1}^{\infty}(Q_1h_j,h_j)_{_{U_1}}<+\infty,$
where $\{h_j\}$ is an orthonormal basis in $U_1$.
\vskip1mm
From the above considerations
$$
\sum_{j=1}^{\infty}(Q_1h_j,h_j)_{_{U_1}}
\le t\sum_{j=1}^{\infty}\|g_j\|_{U_1}^2 \quad \mbox{and then}\quad
\sum_{j=1}^{\infty}(Q_1h_j,h_j)_{_{U_1}}\equiv \mathrm{tr} Q_1<+\infty.
$$
\hfill$\square$
}

\begin{prop} \label{prop2}
For any $a\in U$ the process
\begin{equation}
\left(a,W_c(t)\right)_U=\sum_{j=1}^{\infty}\left(a,g_j\right)_U\beta_j(t)
\end{equation}
is real-valued Wiener process and
$$
E\left(a,W_c(t)\right)_U\left(b,W_c(t)\right)_U=(t\land s)\left(Qa,b\right)_U
\qquad\mbox{for }\ a,b\in U.
$$
Additionally, $\mathrm{Im}\, Q_1^{\frac12}=U_0$ and $\Vert u\Vert_{U_0}=
\left\Vert Q_1^{-\frac12}u\right\Vert_{U_1}$.
\end{prop}

\proof{
We shall prove that the series (3) defining the process
$\left(a,W_c(t)\right)_U$
is convergent in the space $L^2(\Omega,\mathcal{F},P)$.
\vskip2mm
Let us notice that the series (3) is the sum of independent random 
variables with zero mean. Then the series does converge in 
$L^2(\Omega,\mathcal{F},P)$ if and
only if the following series
$\sum_{j=1}^{\infty}E\left(\left(a,g_j\right)_U\beta_j(t)\right)^2
$ converges.
\vskip1mm
Because $J$ is Hilbert--Schmidt operator, we obtain
\begin{eqnarray*}
\sum_{j=1}^{\infty}E\left(\left(a,g_j\right)_U^2\beta_j^2(t)\right)&=&
\sum_{j=1}^{\infty}\left(a,g_j\right)_U^2
\le \|a\|_U^2\sum_{j=1}^{\infty}\left\|g_j\right\|_U^2\\
&\le& C\|a\|_U^2\sum_{j=1}^{\infty}\left\|Jg_j\right\|_{U_1}^2<
+\infty.
\end{eqnarray*}
Hence, the series (3) does converge.
Moreover, when $t\ge s\ge 0$, we have
\begin{eqnarray*}
E\left(\left(a,W_c(t)\right)_U\left(b,W_c(s)\right)_U\right)&=&
E\left(\left(a,W_c(t)-W_c(s)\right)_U\left(b,W_c(s)\right)_U\right)\\
&+&E\left(\left(a,W_c(s)\right)_U\left(b,W_c(s)\right)_U\right)\\
&=&E\left(\left(a,W_c(s)\right)_U\left(b,W_c(s)\right)_U\right)\\
&=&E\left(\left[\sum_{j=1}^{\infty}\left(a,g_j\right)_U\beta_j(s)\right]
\left[\sum_{k=1}^{\infty}\left(b,g_k\right)_U\beta_k(s)\right]\right).
\end{eqnarray*}
Let us introduce
$$
S^a:=\sum_{j=1}^{\infty}(a,g_j)_{_U}\beta_j(t),\quad 
S^b:=\sum_{k=1}^{\infty}(b,g_k)_{_U}\beta_k(t),\ \mbox{for}\ a,b\in U.
$$
Next, let $S_N^a$ and $S_N^b$ denote the partial sums of the series $S^a$
and $S^b$, respectively.
From the above considerations the series $S^a$ and $S^b$ are convergent
in $L^2(\Omega,\mathcal{F},P;\mathbb{R})$.
Hence
$E(S^aS^b)=\lim_{N\to\infty}E(S_N^aS_N^b).$
In fact,
\begin{eqnarray*}
E|S_N^aS_N^b-S^aS^b|&=&E|S_N^aS_N^b-S_n^aS^b+S^bS_n^a-S^bS^a|\\
&\le& E|S_N^a||S_N^b-S^b|+E|S^b||S_N^a-S^a|\\
&\le&\left(E|S_N^a|^2\right)^{\frac12}\left(E|S_N^b-S^b|^2\right)^{\frac12}\\
&+&\left(E|S^b|^2\right)^{\frac12}\left(E|S_N^a-S^a|^2\right)^{\frac12}
~\longrightarrow_{\hspace{-4ex}_{N\to\infty}} 0
\end{eqnarray*}
because $S_N^a$ converges to $S^a$ and $S_N^b$ converges to $S^b$ in
quadratic mean.
\vskip1mm
Additionally,
$E(S_N^aS_N^b)=t\sum_{j=1}^{N}(a,g_j)_{_U}(b,g_j)_{_U}$
and when $N\to +\infty$
$$
E(S^aS^b)=t\sum_{j=1}^{\infty}(a,g_j)_{_U}(b,g_j)_{_U}.
$$
Let us notice that
\begin{eqnarray*}
\left(Q_1a,b\right)_{U_1}&=&E\left(a,W_c(1)\right)_{U_1}\left(b,W_c(1)
\right)_{U_1}
=\sum_{j=1}^{\infty}\left(a,g_j\right)_{U_1}\left(b,g_j\right)_{U_1}\\
&=&\sum_{j=1}^{\infty}\left(a,Jg_j\right)_{U_1}\left(b,Jg_j\right)_{U_1}
=\sum_{j=1}^{\infty}\left(J^*a,g_j\right)_{U_0}\left(J^*b,g_j\right)_{U_0}\\
&=&\left(J^*a,\ \sum_{j=1}^{\infty}\left(J^*b,g_j\right)g_j\right)_{U_0}
=\left(J^*a,J^*b\right)_{U_0}=\left(JJ^*a,b\right)_{U_1}.
\end{eqnarray*}
That gives
$Q_1=JJ^*.$
In particular
\begin{equation}
\left\|Q_1^{\frac12}a\right\|_{U_1}^2=\left(JJ^*a,a\right)_{U_1}
=\left\|J^*a\right\|_{U_0}^2,\qquad a\in U_1.
\end{equation}
\vskip2mm
Having (4), we can use theorems about images of linear operators (e.g.,
\cite{3}, Appendix B.2, Proposition B.1 (ii)).
\vskip1mm
By this theorem
$
\mathrm{Im}\, Q_1^{\frac12}=\mathrm{Im}\, J.$
But for any $j\in\mathbb{N}$, and $g_j\in U_0$, $Jg_j=g_j$, 
that is $\mathrm{Im}\, J=U_0$.
Then $\mathrm{Im}\, Q_1^{\frac12}=U_0.$
\vskip1mm
Moreover, the operator $G=Q_1^{-\frac12}J$ is a bounded operator from $U_0$
on $U_1$. From (4) the joint operator
$G^*=J^*Q_1^{-\frac12}$
is an isometry, so $G$ is isometry, too. Then
$$
\left\|Q_1^{-\frac12}u\right\|_{U_1}=\left\|Q_1^{-\frac12}Ju\right\|_{U_1}=
\|u\|_{_{U_0}}.
$$
\hfill$\square$
}

In the case when $Q$ is nuclear operator, $Q^{\frac12}$ is Hilbert-Schmidt
operator. Taking $U_1=U$, the process $W_c(t)$, $t\ge 0$, defined by (2) is
the classical Wiener process introduced in Definition 1.
\vskip5mm
\begin{defin} \label{def3}
The process $W_c(t)$, $t\ge 0$, defined in (2), is called {\tt cylindrical
Wiener process} in $U$ when $\mathrm{tr} Q=+\infty$.
\end{defin}

The stochastic integral with respect to cylindrical Wiener process is defined
as follows.
\vskip2mm
As we have already written above, the process $W_c(t)$ defined by (2)
is a Wiener process in
the space $U_1$ with the covariance operator $Q_1$ such
that $\mathrm{tr} Q_1<+\infty$. Then the stochastic integral
$\int_0^tg(s)dW_c(s)\in Y$, where $g(s)\in L(U_1,Y)$, with respect to the
Wiener process $W_c(t)$ is well defined on $U_1$.
\vskip2mm
Let us notice that $U_1$ is not uniquely determined. The space $U_1$ can be
an arbitrary Hilbert space such that $U$ is continuously embedded into $U_1$
and the embedding of $U_0$ into $U_1$ is a Hilbert-Schmidt operator. We would
like to define the stochastic integral with respect to cylindrical Wiener
proces $W_c(t)$ (given by (2)) in such a way that the integral is well
defined on the space $U$ and does not depend on the choice of the space
$U_1$.
\vskip5mm
We denote by $N(Y)$ the space of all stochastic processes
\begin{equation}
\Phi:\ [0,T]\times \Omega \to L_2(U_0,Y)
\end{equation}
such that
\begin{equation}
E\left(\int_0^T\left\|\Phi (t)\right\|_{L_2(U_0,Y)}^2dt\right)<+\infty
\end{equation}
and for all $u\in U_0$, $\Phi (t)u$ is a $Y$--valued stochastic process
measurable with respect to the filtration $(\mathcal{F}_t)$.
\vskip2mm
The stochastic integral $\int_0^t\Phi (s)dW_c(s)\in Y$ with respect to
cylindrical Wiener process, given by (2) for any process $\Phi \in N(Y)$,
can be defined as the limit
\begin{equation}
\int_0^t\Phi (s)dW_c(s)=\lim_{m\to \infty}\sum_{j=1}^m
\int_0^t\Phi (s)g_jd\beta_j(s)\qquad \mbox{in}\ Y
\end{equation}\
in $L^2(\Omega)$  sense.

\noindent
{\bf Comment:} Before we prove that the stochastic integral given
by the formula (7) is well defined, let us recall properties of the
operator $Q_1$. From Proposition 1, cylindrical Wiener process $W_c(t)$
given by (2) has the covariance operator $Q_1:\ U_1 \to U_1$, which
is a nuclear operator in the space $U_1$, that is $\mathrm{tr} Q_1<+\infty$.
Next, basing on Proposition 2, $Q_1^{\frac12}:\ U_1 \to U_0$,
$\mathrm{Im}\, Q_1^{\frac12}=U_0$ and 
$\|u\|_{U_0}=\left\|Q_1^{-\frac12}u\right\|_{U_1}$
for $u\in U_0$.

 Moreover, from the above considerations and properties of the operator
$Q_1$  we may deduce that $L(U_1,Y)\subset L_2(U_0,Y)$.
This means that each
operator $\Phi \in L(U_1,Y)$, that is linear and bounded from $U_1$ into $Y$,
is Hilbert-Schmidt operator acting from $U_0$ into $Y$, that is
$\Phi \in L_2(U_0,Y)$ when $\mathrm{tr} Q_1<+\infty$ in $U_1$. This means that
conditions (5) and (6) for the family $N(Y)$ of integrands are natural
assumptions for the stochastic integral given by (7).

Now, we shall prove that the series from the right hand side of (7) is
convergent.
\vskip2mm
Denote
$$
W_c^{(m)}(t):=\sum_{j=1}^mg_j\beta_j(t)
$$
and
$$
Z_m:=\int_0^t\Phi (s)W_c^{(m)}(s),\qquad t\in [0,T].
$$
Then, we have
\begin{eqnarray*}
E\left(\left\|Z_n-Z_m\right\|_Y^2\right)
& = & E\left\|\sum_{j=m+1}^n\int_0^t\Phi (s)g_jd\beta_j
(s)\right\|_Y^2\qquad\mbox{for }\ n\ge m\ge 1\\
&\le & E\sum_{j=m+1}^n\int_0^t\left\|\Phi(s)g_j\right\|_Y^2ds
~~\longrightarrow_{\hspace{-4.5ex}_{m,n\to\infty}}  0,
\end{eqnarray*}
because from the assumption (6)
$$
E\int_0^t\left(\sum_{j=1}^{\infty}\left\|\Phi (s)g_j\right\|_Y^2\right)ds<
+\infty.
$$
\vskip5mm
Then, the sequence $(Z_m)$ is Cauchy sequence in the space of
square--integrable random variables. So, the stochastic integral with respect
to cylindrical Wiener process given by (7) is well defined.
\vskip2mm
As we have already mentioned, the space $U_1$ is not uniquely determined.
Hence, the cylindrical Wiener proces $W_c(t)$ defined by (2)
is not uniquely determined either.
\vskip2mm
Let us notice that the stochastic integral defined by (7) does not depend
on the choice of the space $U_1$. Firstly, in the formula (7) there are not
elements of the space $U_1$ but only $\{g_j\}$--basis of $U_0$.
Additionally, in (7) there are not eigenfunctions of the covariance operator
$Q_1$.
Secondly, the class $N(Y)$ of integrands  does not depend on the choice of
the space $U_1$ because (by Proposition 2.) the spaces $Q_1^{\frac12}(U_1)$
are identical for any spaces $U_1$:
$$
Q_1^{\frac12}:\ U_1\to U_0\qquad\mbox{and}\qquad 
\mathrm{Im}\, Q_1^{\frac12}=U_0.
$$

Hence, the stochastic integral with respect to infinite dimensional Wiener
process, even cylindrical, can be obtained in the above sense as the limit
of stochastic integrals with respect to real-valued Wiener processes.

\section{Connection with Walsh integral}
 
In this section we compare the integral defined in the previous section
with the integral constructed by Walsh \cite{9}.
\vskip2mm
Let us recall from \cite{9} the necessary definitions. Assume that
$(E,\mathcal{E})$ is a Lusin space, i.e. a measurable space homeomorphic
to a Borel subset of the line. (Let us notice that this space includes all
Euclidean spaces and, more generally, all Polish spaces.) Suppose
$\mathcal{A}\subset\mathcal{E}$ is an algebra.

\begin{defin} \label{def4}
Let $(\mathcal{F}_t)$ be a right continuous filtration. A process
$\left\{M_t(A),\mathcal{F}_t,t\ge 1,\right.$
$\left.A\in\mathcal{A}\right\}$
is a {\tt martingale measure} if
\begin{enumerate}
\item $M_0(A)=0$,
\item if $t>0$, $M_t$ is a $\sigma$--finite $L^2$--valued measure,
\item $\left\{M_t(A),\mathcal{F}_t, t\ge 0\right\}$ is a martingale.
\end{enumerate}
\end{defin}

\begin{defin} \label{def5}
A martingale measure $M$ is {\tt othogonal} if, for any two 
disjoint sets $A$ and
$B$ in $\mathcal{A}$, the martingales 
$\left\{M_t(A),\mathcal{F}_t,t\ge 1\right\}$
and $\left\{M_t(B),\mathcal{F}_t,t\ge 1\right\}$ are orthogonal.
\end{defin}

Let us notice that an example of an orthogonal martingale measure is a
white noise. If $W$ is a white noise on $E\times\mathbb{R}_+$, define
$M_t(A)=W\left(A\times[0,t]\right)$. This is clearly martingale measure,
and if $A\cap B=\varnothing$, $M_t(A)$ and $M_t(B)$ are independent,
hence orthogonal.
\vskip2mm
We know how to integrate over $dx$ for fixed $t$ -- this is the Bochner
integral -- and over $dt$ for fixed sets $A$ -- this is the It\^o
integral.
The problem is to integrate over $dx$ and $dt$ at the same
time. Unfortunately, it is not possible to construct a stochastic
integral with respect to all martingale measures. We shall add some
conditions and define a new class of martingale measures.

\begin{defin} \label{def6}
A martingale measure $M$ is {\tt worthy} if there exists a random 
$\sigma$--
finite measure $K(\Lambda,\omega)$, 
$\Lambda\in\mathcal{E}\times\mathcal{E}\times
\mathcal{B}$, where $\mathcal{B}$ is Borel sets on $\mathbb{R}_+$, 
$\omega\in\Omega$,
such that
\begin{enumerate}
\item $K$ is positive definite and symmetric in $x$ and $y$,
\item for fixed $A$, $B$, $\left\{K(A\times B\times(0,t]),
t\ge 0\right\}$ is predictable,
\item for all $n\in\mathbb{N}$, $\mathbb{E}\left\{K\left(
E_n\times E_n\times[0,T]\right)\right\}<+\infty$, where
$E_n\subset\mathcal{E}$,
\item for any rectangle $\Lambda$, $|M(\Lambda)|\le K(\Lambda)$.\\
(We call $K$ the {\tt dominating measure} of $M$.)
\end{enumerate}
\end{defin}

Let us notice, that conditions of the above definition are satisfied by
orthogonal martingale measures, that is orthogonal martingale measures
are worthy.
\vskip2mm
As usual, we first define the integral for elementary functions, then
for simple functions, and then for all functions in a certain class.
\vskip5mm
\begin{defin} \label{def7}
 function $f$ is {\tt elementary} if it is of the form
\begin{equation}
f(s,x,w)=I_{(a,b]}(s)I_A(x)X(w),
\end{equation}
where $0\le a\le t$, $X$ is bounded and  $\mathcal{F}$ -- measurable, 
and $A\in\mathcal{E}$.

A function $f$ is {\tt simple} if it is a finite sum of elementary 
functions.
\end{defin}

We shall denote the class of simple functions by $\mathcal{S}$.

\begin{defin} \label{def8}
The {\tt predictable $\sigma$-field} $\mathcal{P}$ on
 $\Omega\times E\times\mathbb{R}_+$
is the $\sigma$-field generated by~$\mathcal{S}$.
A function is {\tt predictable} if it is $\mathcal{P}$--measurable.
\end{defin}

We define a norm $\|\cdot\|_M$ on the predictable functions by
$$
\|f\|_M=\mathbb{E}\left\{\left(|f|,|f|\right)_K\right\}^{\frac12},
$$
where
$$
(f,g)_K=\int_{E\times E\times\mathbb{R}_+}f(s,x)g(s,y)K(dxdyds).
$$
\vskip2mm
Let $\mathcal{P}_M$ be the class of all predictable $f$ for which
$\|f\|_M<+\infty$.

\begin{prop} \label{prop3}
The class $\mathcal{P}_M$ is a Banach space. Moreover, $\mathcal{S}$ 
is dense in $\mathcal{P}_M$.
\end{prop}

(For proof and details, see \cite{9}.)
\vskip3mm
Now, we can follow Walsh and define stochastic integral as a martingale
measure.
\vskip2mm
If $f$ is an elemntary function, that is $f$ has the form (8), define a
martingale measure $f\cdot M$ by
$$
f\cdot M_t(B)\stackrel{\mbox{df}}{=}
X(w)\left[M_{t\wedge b}(A\cap B)-M_{t\wedge a}(A\cap B)\right].
$$
\vskip5mm
\begin{prop} \label{prop4}
(Lemma 2.4, \cite{9})
The martingale measure $f\cdot M$ is worthy. Moreover
\begin{equation}
\mathbb{E}\left\{\left(f\cdot M_t(B)\right)^2\right\}\le \|f\|_M^2\quad
\mbox{for all $B\in\mathcal{E},\ t\le T$.}
\end{equation}
\end{prop}
 
Now, we can define $f\cdot M$ for $f\in \mathcal{S}$ by linearity.
\vskip2mm
Suppose that $f\in\mathcal{P}_M$. By Proposition 3 there exist 
$f_n\in\mathcal{S}$
such that $\|f-f_n\|_M\to 0$ when $n\to\infty$. By (9), if $A\in\mathcal{E}$
and $t\le T$,
$$
\mathbb{E}\left\{\left(f_m\cdot M_t(A)-f_n\cdot M_t(A)\right)^2\right\}
\le \|f_m-f_n\|_M\to 0,\ \mbox{when $m,n\to\infty$.}
$$
\vskip2mm
It follows that $\left(f_n\cdot M_t(a)\right)$ is Cauchy sequence in
$L^2(\Omega,\mathcal{F},P)$, then it converges in $L^2$ to a martingale
which we shall denote by $f\cdot M_t(A)$. Additionally, the limit is
independent of the choice of the sequence $(f_n)$.
\vskip5mm
\begin{prop} \label{prop5}
(Theorem 2.5, \cite{9}) If $f\in\mathcal{P}_M$, then $f\cdot M$ is a worthy
martingale measure.
\end{prop}

Now, because the stochastic integral is defined as a martingale measure,
we define the "usual" stochastic integrals by
$$
\int_0^t\int_Af(s,x)M(dxds)=f\cdot M_t(A)
$$
and
$$
\int_0^t\int_Ef(s,x)M(dxds)=f\cdot M_t(E).
$$
\vskip3mm
Let us consider the integral constructed by Walsh and recalled in
this section in the case when $E\equiv \mathbb{R}^d$ and $M_t(A)$ is
cylindrical Wiener process. In this case
\vskip2mm
$$
M_t(A)\equiv W(t,A,w),
$$
where $t\ge 0$, $A\in\mathcal{B}(\mathbb{R}^d)$ with $mes(A)<+\infty$, and
$w\in\Omega$. (In the remaining part of the paper we shall omit the
argument $w$.)
\vskip2mm
For any $A$ fixed, $\{W(\cdot,A)\}$ is a real--valued Wiener process,
adapted to the filtration $(\mathcal{F}_t)$ which does not depend on A.
Moreover,
$\mathbb{E}\left(W(t,A)W(s,B)\right)=t\wedge s\ mes(A\cap B)$.
\vskip2mm
Assume that
$\phi:\ [0,\infty)\times\mathbb{R}^d\times\Omega\to\mathbb{R}.$
\vskip2mm
Using Walsh approach we can define the integral
$$
J(\phi):=\int_0^T\int_{\mathbb{R}^d}\phi(t,x)W(dt,dx).
$$
\vskip2mm
We start from simple functions of the form
$\phi(t,x,w)=f(t,w)I_A(x)$.
Then we have
$$
J(\phi)=\int_0^Tf(t)W(dt,A)\equiv\int_0^Tf(t)dW(t,A).
$$
\vskip2mm
Let us introduce the following classes of functions.
\vskip2mm
By
$P_T$ we denote the class of functions
$\phi: [0,T]\times\mathbb{R}^d\times\Omega\to\mathbb{R}$
satisfying the following conditions:
\begin{enumerate}
\item $\phi$ is measurable,
\item the function $\phi(t,x)$ is
$\mathcal{F}_t$--measurable,
\item $\mathbb{E}\left(\int_0^T\int_{\mathbb{R}^d}|\phi(t,x)|^2dtdx
\right)<+\infty.$
\end{enumerate}
By $\widetilde P_T$ we denote the class of functions
$\phi:\ [0,T]\times\Omega\to L^2(\mathbb{R}^d)$
such that:
\begin{enumerate}
\item $\phi$ is measurable,
\item the function $\phi(t)$ is
$\mathcal{F}_t$--measurable,
\item $\mathbb{E}\left(\int_0^T|\phi(t)|^2_{L^2(\mathbb{R}^d)}dt
\right)<+\infty.$
\end{enumerate}
Let us notice that the both classes $P_T$ and $\widetilde P_T$ coincide.
\vskip2mm
Now, we can formulate the following result.
\begin{prop} \label{prop6}
Let $W(t)$ be a cylindrical Wiener process. Then
\begin{equation}
\int_0^T\phi(t)dW(t)=\int_0^T\int_{\mathbb{R}^d}\phi(t)(x)W(dt,dx)
\end{equation}
for any $\phi\in\widetilde P_T$.
\end{prop}
\proof{
It is enough to check the formula (10) for simple function
$\phi(t)(x)\equiv f(t)I_A(x)$. Let $\{e_k\}$, $k=1,2,...,d$, be
a basis in $\mathbb{R}^d$, where  $e_k=I_A/mes(A)$ and
$e_k$, for $k=2,...,d$, are arbitrary.
\begin{eqnarray*}
\int_0^T\int_{\mathbb{R}^d}\phi(t)(x)W(dt,dx)
&\equiv&\int_0^T\int_{\mathbb{R}^d}f(t)I_A(x)W(dt,dx)\\
&\equiv&\sum_{k=1}^{d}\int_0^T(f(t),e_k)W(dt, e_k)\\
&\equiv&\sum_{k=1}^{d}\int_0^T\left(f(t),e_k\right)dW(t)[e_k]
=\int_0^T\phi(t)dW(t).
\end{eqnarray*}
\hfill $\square$
}

\noindent{\bf Comment:} Another, very recent example of
integrating over random measures in
multidimensional spaces has been given in the paper of Peszat
and Zabczyk \cite{8}. In the paper, the noise is supposed
to be a spatially homogeneous Wiener
process in
some special space. The authors describe its reproducing kernel
and provide
the concept of stochastic integral with respect to introduced Wiener
process.

{\bf Acknowledgements.} I would like to thank  M. Capi\'nski
for helpful comments and suggestions.


\begin{thebibliography}{AB}
\bibitem{1} Balakrishnan A.\
\textit{Applied Functional Analysis}, Springer--Verlag, New York, 1981.

\bibitem{2} Capi\'nski M., Peszat S.\ \textit{
Local existence and uniqueness of strong solutions to 3--D stochastic
Navier--Stokes equations} Nonlin.\ Diff.\ Eqns.\ Appl., 
{\bf 4}, No. 2, (1997), 185--200.

\bibitem{3} DaPrato G., Zabczyk J.\ \textit{
Stochastic equations in infinite dimensions},
Cambridge University Press, Cambridge, 1992

\bibitem{4} Flandoli F., G\c atarek D.\ \textit{
Martingale and stationary solutions for stochastic
Navier--Stokes equations},
Probab.\ Theory Relat.\ Fields, {\bf 102}, No 3 (1995), 367--391.

\bibitem{5} Gelfand I., Vilenkin N.\  \textit{
Some Problems of Harmonic Analysis, (in Russian)},
Fizmatgiz, Moskwa, 1961.

\bibitem{6} Ichikawa A.\  \textit{
Stability of semilinear stochastic evolution equations},
J.\ Math.\ Anal.\ Appl.,{\bf 90} (1982), 12-44.

\bibitem{7} Karczewska A.\  \textit{
Stochastic solutions to turbulent diffusion},
Nonlinear Analysis, {\bf 37} (1999), 635-675.

\bibitem{8} 
Peszat S., Zabczyk J.\  \textit{
Stochastic evolution equations with a spatially homogeneous
Wiener process},  Stochastic Process.\ Appl.  {\bf 72}, No. 2, (1997), 
 187--204.


\bibitem{9}
Walsh J.\  \textit{
An introduction to stochastic partial differential
equations},
Ecole d'Et\'e de Prob. de St--Flour XIV--1984,
Lecture Notes Math., Vol. 1180,
Springer--Verlag, Berlin--New York, 1986.

\end{thebibliography}
\end{document}